\documentclass{amsart}
\usepackage{amsmath, amssymb}
\usepackage{float, mathtools, listings, amsthm}
\usepackage[f]{esvect}
\usepackage[margin=1in]{geometry}
\usepackage{algpseudocode}
\usepackage{float, mathtools, listings, amsthm}
\usepackage[f]{esvect}
\usepackage[margin=1in]{geometry}
\usepackage{algorithmicx}
\usepackage{algorithm}
\usepackage{algpseudocode}
\usepackage{wrapfig}
\usepackage[numbers]{natbib}

\title{Use of Statistically Leinert Sets to Calculate Return Probabilities of Random Walks in $\mathbb{F}_{s_1} \times \mathbb{F}_{s_2}$}
\author{Colton Griffin, Sanchita Chakraborty}
\date{June 2021}

\newtheorem{thm}{Theorem}[section]

\newtheorem*{thm*}{Theorem}

\newtheorem{prop}[thm]{Proposition}
\newtheorem*{prop*}{Proposition}

\theoremstyle{definition}
\newtheorem{defn}[thm]{Definition}

\usepackage{natbib}
\usepackage{graphicx}

\begin{document}
\maketitle
\section{Abstract}
Hastings first presented bounds on the second largest eigenvalue for matrices in a Hermitian complete positive map in 2007. In this work we extend his work to tighten these bounds. To do this, we introduce the idea of \textit{Statistically Leinert Sets} to modify the generating functions presented in Woess in 1986 and recompute the radii of convergence in his paper in 1986. We primarily use techniques from combinatorics and calculate norms using the ideas presented by Akemann and Ostrang in their paper in 1976.
\section{Introduction}
\subsection{Literature Review}
 \begin{wrapfigure}{r}{0.25\textwidth}
\begin{center}
    \includegraphics[width=0.25\textwidth]{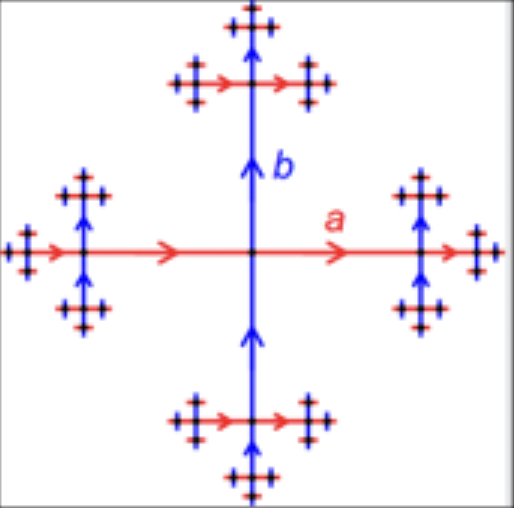}
    \caption{Random Walk with 2 Generators}
\end{center}
\end{wrapfigure}

Random walks have been studied extensively, but some interesting properties arise when restricting to trees. A tree in graph theory is defined as a finite graph which is undirected and connected, lacks cycles, and stems from a root. As we are considering an unbiased walk, so in terms of the elements of the string, they are equally likely to appear. In the graph sense, the movement in any direction is equally likely. In Figure 1, a sample random walk generated by {a, b} is shown. Random walks in $\mathbb{F}_{s_1} \times \mathbb{F}_{s_2}$ which we consider here in this paper can be thought of as a product of one dimension walks. As dimensions increase, the walk can become more complex. An intuitive question then arises: what is the return probability of the walk to the origin? It has been established that random paths on mathematical trees with a number n generators can be expressed using matrices. These paths can either return to the origin or diverge with no way back. This divergence behavior can be understood by studying the radius of convergence. 

In 1976, Akemann and Ostrang \cite{akemann} introduced norms of operators on the Hilbert space: $L^2(G)$, where elements can be written as generally infinite linear combinations acting on the space by left multiplication. In 1986, Woess \cite{woess} presented work on calculating Norms of Free Convolution Operators, using generating functions to represent return probabilities of random walks on a homogeneous tree of degree 2s. In 2007, M.B. Hastings \cite{hastings} first presented bounds on Haar Random Matrices. His work analyzed random walks on a Cayley Tree and specifically considered it for the Hermitian Case. A Cayley Tree \cite{cayley} is an interesting case of trees in graphs, where each non-leaf vertex (i.e. a vertex with no children vertices) has a constant number of vertices stemming from it. For example, a 2-Cayley tree will have 2 vertices stemming from each leaf vertex. This work by Hastings provided the main motivation for the paper presented here.

\subsection{Approach}
For the purposes of this paper the walks are unbiased. In this paper these paths were represented using ”words”, i.e. strings generated from an n number of generators. The set of all possible words can be represented as a set S which forms an algebraic group: the free group: F. We modified a key property, known as the Leinert property, to be imposed on this group. We instead applied a Statistically Leinert Property.

\section{Product of Free Groups}

\begin{defn}
Given a subset $X$ of a group $G$, a string $S$ of elements (or inverses of elements) from $X$ is \textit{valid} if it is of the form
\[S=x_{1}^{-1}x_2\ldots x_{2n-1}^{-1}x_{2n},\]
where $x_i\in X$ for all $i$ and $x_{i}\neq x_{i+1}$.

A string $S$ is \textit{reduced} if it is of the form ($x_i\in X$ and $\varepsilon_i\in \{-1,1\}$)
\[S=x_1^{\varepsilon_1}x_2^{\varepsilon_2}\ldots x_k^{\varepsilon_k},\]
and for $1\leq i<k$ we have either $x_i\neq x_{i+1}$ or $\varepsilon_i=\varepsilon_{i+1}$.

A string $S$ is \textit{bad} if it is valid or reduced and reduces to $S=e$, where $e$ is the identity element in $G$.
\end{defn}

\begin{defn}
A subset $X\subset G$ is a \textit{statistically Leinert set} if the probability of a valid string $S$ of length $2n$ being bad approaches 0 as $n\to\infty$. An equivalent condition is that if $b_{2n}$ is the number of bad valid strings of length $2n$, then
\[\lim_{n\to\infty}\frac{b_{2n}}{s(s-1)^{2n-1}}=0,\]
where $s=|X|$. In particular the growth rate of bad strings is: $lim_{n \rightarrow \infty} sup (\frac{b_{2n}}{s(s-1)^{2n-1}})^{1/n} = 1$
\end{defn}
This equivalent condition comes from the fact that the total number of valid strings is $s(s-1)^{2n-1}$, so the frequency of bad strings is $b_{2n}$ divided by this.
\begin{defn}
A subset $X\subset G$ is a \textit{statistically free set} if the probability of a reduced string $S$ of length $k$ being bad approaches 0 as $k\to\infty$. An equivalent condition is that if $b_n$ is the number of bad reduced strings of length $n$, then 
\[\lim_{n\to\infty}\frac{b_n}{2s(2s-1)^{n-1}}=0,\]
\end{defn}

A Leinert set is statistically Leinert as well. The same follows for free sets.
\begin{prop}
A subset $X\subset G$ is statistically free iff $e\notin X$ and $X\cup \{e\}$ is statistically Leinert.
\end{prop}
\begin{proof}
Suppose $X$ is statistically free. If $e\in X$, then we can consider any string 
\end{proof}
Now a question: How do we show that $b_{2n}$ satisfies the limit above? If we can show that it has an exponential growth rate slower than $(s-1)^2$ (which is the exponential growth rate of the total number of valid strings), then we are done.

\section{Counting Bad Strings}
\begin{defn}
If $S$ is a valid string of the form
\[S=x_{1}^{-1}x_2\ldots x_{2n-1}^{-1}x_{2n},\]
then define the \textit{conjugate} string $\bar S$ to be
\[\tilde S=x_{1}x_2^{-1}\ldots x_{2n-1}x_{2n}^{-1}.\]
In general, if $S$ is any string of letters, then the conjugate string is given by replacing each letter with its inverse.
\end{defn}
We will think about how to count bad strings by considering substrings of larger strings:
\begin{defn}
A \textit{substring} $R\leq S$ is a subset of the letters of the string $S$ given by removing all the letters to the left of the lowest index letter and all the letters to the right of the highest index letter. A bad substring $R\leq S$ is a \textit{kernel substring} if it has no bad substrings.
\end{defn}

Before talking about bad strings, it is worth considering a couple counting problems. For example, what is the number of different random walks on the free group $\mathbf{F}_s$ such that we start and end at the identity? First let's count the different ways to add natural numbers up to $N$, considering ordering as well. This is different to the partition number, which doesn't care about ordering. Note that the total number of ways to satisfy the identity $\ell_1+\ldots+\ell_k=n$, where $\ell_i\geq0$, is $\binom{n+k-1}{k-1}$. If we add $k$ to each side and let $l_i=\ell_i+1$, then we get $l_1+\ldots+l_k=n+k=N$, which is the sum of natural numbers to equal $N$ is $\binom{N-1}{k-1}$. Now we sum over each combination, and we get
\[\sum_k\sum_{\substack{l_1+\ldots+l_k=N\\\ell_i\geq 1}}=\sum_{k=1}^{N}\binom{N-1}{k-1} = 2^{N-1}.\]
Now we count the number of random walks similarly; a walk of length $2n$ can be partitioned into mini-walks back to the identity, where the lengths of all the walks add up to $2n$. Each walk has even length, so we partition the walk length into $2l_1+\ldots+2l_k=2n$. The number of different walks from the identity to the identity of length $2l$ where the walk never touches the identity is $2s(2s-1)^{l-1}$, so we have
\[\text{\# of walks} = \sum_{k=1}^{N}\sum_{\substack{l_1+\ldots+l_k=N\\\ell_i\geq 1}}\prod_{i=1}^{k}2s(2s-1)^{l_i-1} = \sum_{k=1}^{N}\sum_{\substack{l_1+\ldots+l_k=N\\\ell_i\geq 1}}\left(\frac{2s}{2s-1}\right)^k(2s-1)^{N}\]
\[=(2s-1)^{N}\sum_{k=1}^{N}\binom{N-1}{k-1}\left(\frac{2s}{2s-1}\right)^k=(2s-1)^N\frac{2s}{4s-1}\left(\frac{4s-1}{2s-1}\right)^N=2s(4s-1)^{N-1}.\]
Similarly, we can sum over the number of times we return to the identity for the bad strings.

\subsection{Lattice Group}
Here we consider $G=\mathbf{F}_1\times\mathbf{F}_1\times\mathbf{F}_1=\mathbb{Z}^3$, with $X=\{x,y,z\}$ the set of generators for the group. Note that the smallest bad string has length 6, with an example given by
\[S=x^{-1}yz^{-1}xy^{-1}z.\]
\begin{itemize}
    \item For $2n=6$, we simply consider every possible arrangement of strings of the above kind. This is just the total number of permutations of a set of three elements, so $b_6=3!=6$ total bad strings.
    \item For $2n=8$, we cannot have kernel substrings of length 8, so the only possibility is to conjugate by an element to get a string $S_{10}=w^{-1}\tilde S_8w$. $w$ is determined by the ends of $S_8$, leaving only one option per string. So again, there are exactly $b_8=6$ bad strings of length 8.
    \item For $2n=10$, we expect that there are 42 bad strings in total. We can construct these in a few ways. First, we take any bad string of length 8 and append a letter to both sides, yielding 12 options. The rest are given by considering 
\end{itemize}
Here is an attempt to create an upper bound. Every string is represented with three walks on each axis of the lattice. For the positive steps, we have a total of $n=n_x+n_y+n_z$ steps to travel. We pick where these go, giving a $\binom{n}{n_x,n_y,n_z}$ factor. For the remaining spots, we have to decide where the other factors go. 
\subsection{Product of Groups}
Here we consider $G=\mathbf{F}_{s_1}\times \mathbf{F}_{s_2}$, with $X=\{x_1,\ldots,\}$. Note that the smallest bad string has length 8, with an example given by
\[S=x_1^{-1}x_2y_1^{-1}y_2x_2^{-1}x_1y_2^{-1}y_1.\]
\begin{itemize}
    \item For $2n=8$, the total number of bad strings is $\boxed{b_{8}=2s_1(s_1-1)s_2(s_2-1),}$ as we can pick anything for $x_1$ and $y_1$, with $x_2$ and $y_2$ constrained such that they are not $x_1$ and $y_1$ respectively. We can also invert the string, providing the factor of 2.
    \item For $2n=10$, we can use every single string from before and append elements to the end that cancel each other out. Using the kernel form for $S_8$ (string of length 8), we can let $\bar S = z_1^{-1}S_{8}z_1$, where $z_1\in X$ and $z_1\neq x_1$ or $y_2$. So we have $s-2$ choices for $z_1$. There are no other strings of this size that we can consider.
    \item A general strategy for counting the total number of bad strings of a given length is to sum over the size of the smallest bad substring. 
\end{itemize}

\section{Woess' argument}
Here are some attempts to write down similar equations for Woess' bounds. The original set of equations are
\begin{align*}
    \mu^{(n)}(e)&=\sum_{k=0}^{n}f^{(2k)}\mu^{(n-k)}(e)+\alpha_0p^{(2n-1)},\\
    p^{(2n-1)}&=\sum_{k=0}^{n-1}f^{(2k)}p^{(2n-2k-1)}+\alpha_0\mu^{(n-1)}(e),\\
    a^{(2n)}_i&=\sum_{k=0}^{n}\left(f^{(2k)}-f^{(2k)}_i\right)a^{(2n-2k)}_i + \alpha_0b^{(2n-1)}_i,\\
    b^{(2n-1)}_i&=\sum_{k=0}^{n-1}f^{(2k)}b^{(2n-2k-1)}_i+\alpha_0a^{(2n-2)}_i,\\
    f^{(2n)}_i&=\alpha_i^2a^{(2n-2)}_i.
\end{align*}
Almost all of these identities still hold if we consider the group $G=\mathbf{F}_{s_1}\times\ldots\times\mathbf{F}_{s_N}$. We denote the set of generators $Y$ for $G$ as 
\[Y=\{x_{i,j} : i=1,\ldots,N;\, j=1,\ldots,s_{i}\}.\]
Additionally, we let $X=Y\cup \{e\}$, where $e$ is the identity for $G$.
We write the following definitions:
\begin{align*}
    \mu^{(n)}(e) &= \operatorname{Pr}[X_{2n} = e : X_0 = e],\\
    p^{(2n-1)} &= \operatorname{Pr}[X_{2n} = e : X_1 = e],\\
    f_{i,j}^{(n)} &= \operatorname{Pr}[X_n=e ;\, X_m\neq e\text{ for }m=1,\ldots n-1;\, X_1=x_{i,j}^{-1} : X_0=e],\\
    a_{i,j}^{(2n)} &= \operatorname{Pr}[X_{2n+1}=e ;\, X_m\neq x_{i,j}\text{ for }m=2,\ldots 2n : X_1=e],\\
    b_{i,j}^{(2n-1)} &= \operatorname{Pr}[X_{2n-1}=e ; X_m\neq x_{i,j}\text{ for }m=1,\ldots 2n-2 : X_0=e],\\
    d_{i,j}^{(n)} &= \operatorname{Pr}[X_n=e ;\, X_{n-1}\neq x_{i,j}^{-1};\, X_m\neq e\text{ for }m=1,\ldots n-1;\, X_1=x_{i,j}^{-1} : X_0=e].
\end{align*}
We also write the generating functions
\begin{align*}
G(z) = \sum_{n=0}^{\infty}\mu^{(n)}(e)z^{2n}, \quad H(z)=\sum_{n=1}^{\infty}p^{(2n-1)}z^{2n-1},\quad F_{i,j}(z)=\sum_{n=0}^{\infty}f_{i,j}^{(2n)}z^{2n},\\
    A_{i,j}(z)=\sum_{n=0}^{\infty}a_{i,j}^{(2n)}z^{2n},\quad B_{i,j}(z)=\sum_{n=1}^{\infty}b_{i,j}^{(2n-1)}z^{2n-1},\quad D_{i,j}(z)=\sum_{n=0}^{\infty}d_{i,j}^{(2n)} z^{2n},
\end{align*}
where we write $f^{(n)}=\sum_{i,j}f_{i,j}^{(n)}$ and $F(z)=\sum_{i,j}F_{i,j}(z)$. Once again, $f^{(n)}_{j,i}=0$ if $n$ is odd, which means the conclusion presented in the Woess paper still holds. The first four identities may be interpreted as
\begin{align*}
    \mu^{(n)}(e)&=\sum_{k=0}^{n}f^{(2k)}\mu^{(n-k)}(e)+\alpha_0p^{(2n-1)},\\
    p^{(2n-1)}&=\sum_{k=0}^{n-1}f^{(2k)}p^{(2n-2k-1)}+\alpha_0\mu^{(n-1)}(e),\\
    a^{(2n)}_{i,j}&=\sum_{k=0}^{n}\left(f^{(2k)}-f^{(2k)}_{i,j}\right)a^{(2n-2k)}_{i,j} + \alpha_0b^{(2n-1)}_{i,j},\\
    b^{(2n-1)}_i&=\sum_{k=0}^{n-1}f^{(2k)}b^{(2n-2k-1)}_{i,j}+\alpha_0a^{(2n-2)}_{i,j}.
\end{align*}
For each of these identities, we sum over the index of first return to the identity. The last identity is a bit harder to work with since $d_{i,j}^{(n)}$ is nonzero if $X$ is not Leinert:
\[f_{i,j}^{(2n)} = \alpha_{i,j}^2a_{i,j}^{(2n-2)} + d_{i,j}^{(2n)}.\]
This implies that
\[F_{i,j}(z)=\alpha_{i,j}^2z^2A_{i,j}(z)+D_{i,j}(z).\]
Using the results of Woess, we get the relations
\[G(z)=\frac{1}{1-F(z)-F_0(z)},\quad F_{i,j}(z) = \frac{\alpha_{i,j}^2z^2}{F_{i,j}(z) + 1/G(z)} + D_{i,j}(z)\]
\[\implies F_{i,j}(z) = \frac{\sqrt{(1+D_{i,j}(z)G(z))^2 + 4\alpha_{i,j}^2z^2G(z)^2} +D_{i,j}(z)G(z)-1}{2G(z)}.\]
Therefore we have
\begin{align*}
    G(z) &= 1+\frac{1}{2}\left(\sqrt{1+4\alpha_{0}^2z^2G(z)^2}-1\right)\\
    &+\frac{1}{2}\sum_{i,j}\left(\sqrt{(1+D_{i,j}(z)G(z))^2 + 4\alpha_{i,j}^2z^2G(z)^2} +D_{i,j}(z)G(z)-1\right).
\end{align*}
$$Q\left(t,z\right)=1+\frac{1}{2}\sum_{i,j}\left(\sqrt{\left(1+\frac{zt^{2}}{c^{2}-t^{2}}\right)^{2}+4\alpha_{i,j}^{2}z^{2}t^{2}}+\frac{zt^{2}}{c^{2}-t^{2}}-1\right)$$
Here, $D_{i,j}$ is a power series with coefficients vanishing as $n\to\infty$. More specifically, suppose that $d_{i,j}^{(2n)}\sim \beta^{-n}$ for large $n$. Then we have $R^{-1}=\lim_{n\to\infty}(d_{i,j}^{(2n)})^{1/2n}=\beta^{-1/2}$, and so the radius of convergence is root of the decay rate $\sqrt\beta$ of the frequency of bad strings, where $\beta<(2s-1)^2$.

\subsection{Special Case (Cubic Equation Derivation)}
Consider the case where $\alpha_{ij}=a$ are all the same, say $a=\frac{1}{2s}$, for the group $G=\mathbf{F}_s\times\mathbf{F}_s$. In addition, by symmetry we can argue that the functions $D_{ij}=D$ are all the same. Then we have that $G(z)= Q_R(z,G(z))$ is given by the equation
\[G(z) = 1+s\left(\sqrt{(1+D(z)G(z))^2 + 4a^2z^2G(z)^2} +D(z)G(z)-1\right),\]
where $Q_R$ is to denote replacing $D(z)$ with $\frac{z^2}{R^2-z^2}$, since $D(z)\leq\frac{z^2}{R^2-z^2}$. Therefore $Q(z,t)\leq Q_R(z,t)$.
Moving terms to the other side and squaring gives us
\[\iff [(G-1)-s(DG-1)]^2=(G-1)^2+s^2(DG-1)^2 - 2s(DG-1)(G-1)=s^2[(1+DG)^2 + 4a^2z^2G^2]\]
\[\iff (G-1)^2 - 2s(DG-1)(G-1)-s^2[4DG + 4a^2z^2G^2]=0.\]
This yields a quadratic equation in $G$, namely $AG^2+BG+C=0$, where
\[A = 1 - 2sD - 4a^2z^2s^2,\quad B = -4s^2D+2sD+2s-2,\quad C = 1-2s.\]
Therefore $G=\frac{1}{2A}(-B\pm\sqrt{B^2-4AC})$, and as a power series $G$ returns real values for all values of $z$, so we can solve for the radius of convergence as $B^2-4AC=0$:
\[\left(4s^2D-2sD-2s+2\right)^2 - 4\left(4a^2z^2s^2+ 2sD -1\right)\left(2s-1\right)=0\]
\begin{equation}
\iff \left(4s^2z^2-2sz^2-\left(2s-2\right)\left(R^2-z^2\right)\right)^2 - 4\left(2sz^2+\left(4a^2z^2s^2 -1\right)\left(R^2-z^2\right)\right)\left(R^2-z^2\right)\left(2s-1\right)=0.
\end{equation}
Expanded out, the equation is
\begin{equation}
    - 32 a^2 s^3 z^6 + 16 a^2 s^2 z^6  + 64 a^2 R^2 s^3 z^4 - 32 a^2 R^2 s^2 z^4 + 16 s^4 z^4 - 32 a^2 R^4 s^3 z^2 + 16 a^2 R^4 s^2 z^2 - 16 R^2 s^3 z^2 + 4 R^4 s^2=0.
\end{equation}
This equation in $z$ yields the radius of convergence for $G$, assuming that $z<R$. Here is an example plot comparing the curves $y=\left(2\sqrt{2x-1}\right)^{-1}$ and $(x,y)=(s,z)$, where $(x,y)$ satisfy the cubic equation above:
\begin{figure}[H]
    \centering
    \includegraphics[width=0.8\linewidth]{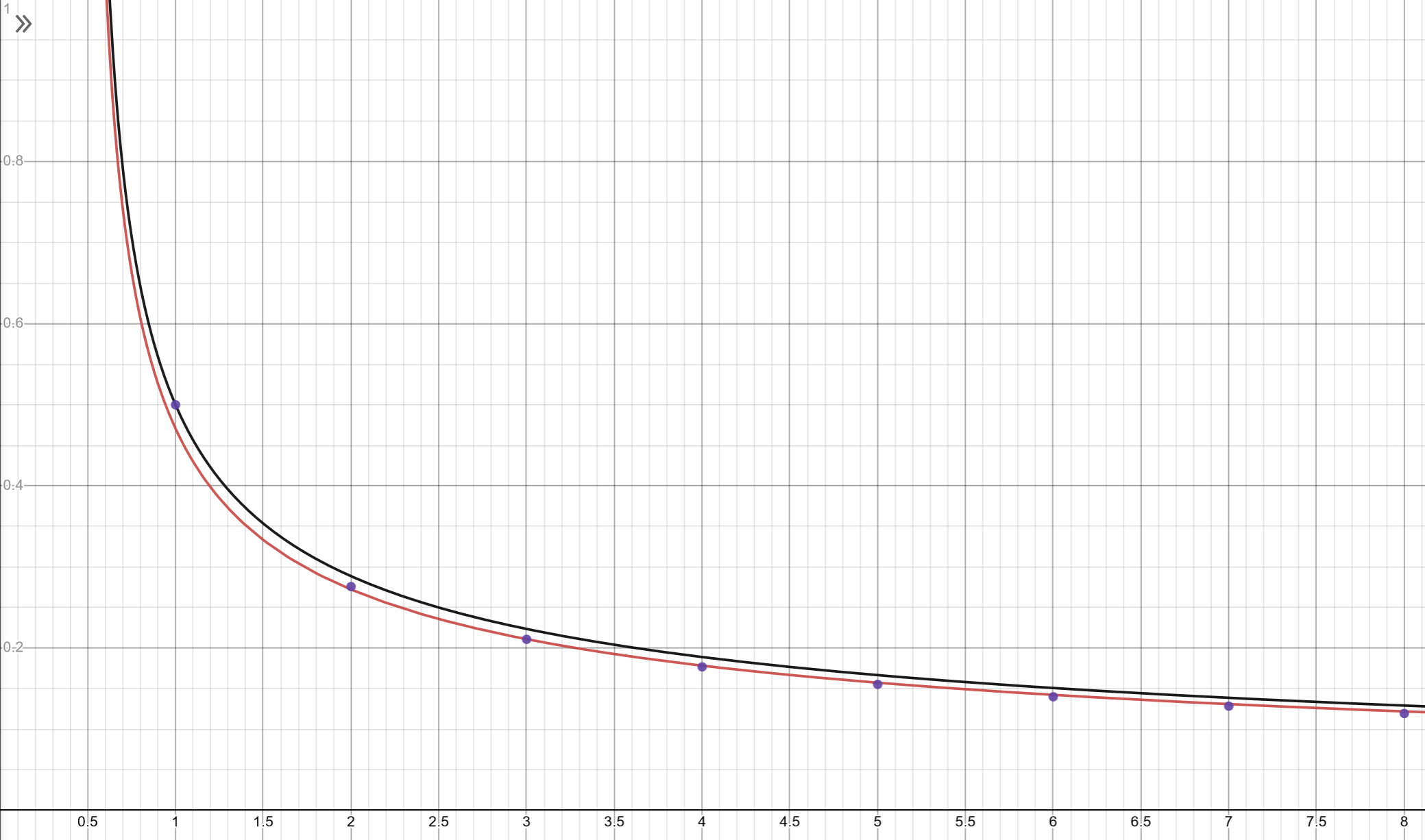}
    \label{fig:my_label}
\end{figure}
The purple points in the plot show numerically computed radii of convergence using the following method. Let $(U_1,\ldots,U_s)\in\operatorname{U}(N)$ and $(V_1,\ldots,V_s)\in \operatorname{U}(N)$ be Haar randomly sampled unitary matrices. Then the dots are given by computing the quantity $z^{-1} = \|a\sum_i U_i\otimes I + I\otimes V_i\|$, with $\|.\|$ the 2-norm. Here each point is given by averaging over 4 computations, $a=1$, and $N=75$.

More generally, if we had some unknown function $D$, then we can use the quadratic formula once more to find the relation for $D$:
\[\left(4s^2D-2sD-2s+2\right)^2 - 4\left(4a^2z^2s^2+ 2sD -1\right)\left(2s-1\right)=0\]
\[\iff 4 a^2 (1 - 2 s) z^2 + ((2 s - 1) D - 1)^2 = 0\iff D(z)=\frac{1\pm 2az\sqrt{2 s - 1 }}{(2s-1)}\]
\[\iff z = \frac{\pm 1}{2a\sqrt{2s-1}}\left(1-(2s-1)D(z)\right).\]
The sign is indeterminate as $D$ cannot be determined. If $D=0$, then we get $z^{-1} = 2a\sqrt{2s-1}$ as expected.

Note that in the above equation, it is quadratic in $R^2$, which we can solve as
\[R^2 = \frac{2\sqrt{a^2\left(2s-1\right)^3z^6} +4 a^2 (2 s - 1) z^4 - 2 s z^2}{4a^2\left(2s-1\right)z^2-1} = \frac{2az^3\left(2s-1\right)^{3/2} +4 a^2 (2 s - 1) z^4 - 2 s z^2}{4a^2\left(2s-1\right)z^2-1}\]

\subsection{Proof for product of free groups}
Let $\beta^{(2n)}$ denote the number of bad strings of length $2n$ for a set $X\subset G$, and let $\mathcal{B}^{(2n)}=\frac{\beta^{(2n)}}{s(s-1)^{2n}}$ be the frequency of bad strings of length $2n$. Say that $G=\mathbf{F}_s\times \mathbf{F}_s$, with $X$ the generator set. We have an upper bound on $\beta^{(2n)}$ that goes to 0 as $s\to\infty$, therefore so does $\beta^{(2n)}$. Recall from the definition of $d_{ij}^{(2n)}$ that only bad strings contribute to a nonzero quantity, so we can say that $D_{ij}(z) = \frac{c^2z^2}{1-c^2z^2}$, where $c$ is the asymptotic exponential rate of encountering these strings. Therefore $c < c_{\text{max}}$, where $c_{\text{max}}$ is the exponential rate of encountering any kind of bad string. Since $c_{\text{max}}\to 0$, so does $c$, and therefore $D_{ij}(z)\to 0$ as we increase $s$. More conveniently, given some fixed uniform constant for all of the coefficients of $\alpha$, we expect that the radius of convergence for $G(z)$ will be $r^{-1} \sim 2a\sqrt{2s-1}$.

Take for example $G=\mathbb{Z}^s$. The norm for $\alpha$ with all coefficients set to 1 is $\|\alpha\| = s$. I expect that the frequency of bad strings does not vanish as we increase $s$, and in fact we have some nice asymptotic relation. 
\subsection{Small Idea about Radius of Convergence}
In the Woess paper, they show that $r^{-1}=P(\theta)/\theta=\min\{P(t)/t : t\geq 0\}$, or $r = \max\{t/P(t) : t\geq 0\}$, where
\[P(t) = 1+\frac{1}{2}\sum_{x\in X}\left(\sqrt{1+4|\alpha(x)|^2t^2}-1\right),\quad P'(t)=\frac{1}{2}\sum_{x\in X}\left[\sqrt{1+4|\alpha(x)|^2t^2}-\frac{1}{\sqrt{1+4|\alpha(x)|^2t^2}}\right]\]
Another approach is to state that since $G(z)=P(zG(z))$, the radius of convergence is bounded by $G(z)$ having a finite derivative, and so if the graph of $G(z)$ has a vertical slope for a finite $G(z)$ value, then that could also be the radius of convergence. Suppose that $y=P(xy)$, then the derivative is given by
\[y' = (y+xy')P'(xy) = (y+xy')\frac{1}{2}\sum_{x\in X}\frac{4|\alpha(x)|^2x^2y^2}{\sqrt{1+4|\alpha(x)|^2x^2y^2}}.\]
Letting $y'\to\infty$, we must have that
\[\frac{y'}{y+xy'} \to \frac{1}{x}=\frac{1}{2}\sum_{x\in X}\left[\sqrt{1+4|\alpha(x)|^2x^2y^2}-\frac{1}{\sqrt{1+4|\alpha(x)|^2x^2y^2}}\right].\]
Note that the first term simplifies to
\[\frac{y'}{y+xy'} \to \frac{1}{x}=y-1+\frac{1}{2}\sum_{x\in X}\left[1-\frac{1}{\sqrt{1+4|\alpha(x)|^2x^2y^2}}\right].\]

One option here is to just explicitly solve for the curve for $G(z)$ where we can. If all the constants are $\alpha_{i,j}$ and we do not have the identity, then this simplifies a lot:
\[G(z)=Q(z,G(z))=1 + \frac{1}{2}\sum_{i,j}(\sqrt{(1+\frac{G(z)z^2}{c^2-z^2})^2+4\alpha_{i,j}^2G(z)^2z^2+\frac{G(z)z^2}{c^2-z^2}-1)}\]

\subsection{Existence of bad strings}
Here are some examples of bad strings, to show that they exist for groups of the form $G=\mathbf{F}_{s_1}\times\ldots\times\mathbf{F}_{s_m}$, where $s_i\in \mathbb N$ for all $i$. If $m=1$, then the set $X$ of generators for $G$ is a Leinert set, as they have no nontrivial algebraic relations. If $m\geq 2$, then we run into problems.
\begin{enumerate}
     \item Consider $G=\mathbf{F}_{s_1}\times\mathbf{F}_{s_2}$, which has the generator set $X=\{x_1,\ldots,x_{s_1},y_1,\ldots,y_{s_2}\}$. Here in this group the $x$'s commute with the $y$'s. If $s_1$ and $s_2$ are greater than 1, then we can let $S$ be the string
     \[S=x_1^{-1}x_2y_1^{-1}y_2x_2^{-1}x_1y_2^{-1}y_1=e.\]
     \item Consider $m\geq 3$, then for this group we can consider the string
     \[S=x_1^{-1}y_1z_1^{-1}x_1y_1^{-1}z_1=e,\]
     where $x_1$, $y_1$, and $z_1$ are generators associated with three different copies of a free group in the total group $G$.
 \end{enumerate}
 There are some cases that do provide Leinert sets.
 \begin{enumerate}
    \item If $G=\mathbf{F}_1\times\mathbf{F}_1$, then the generator set $X=\{x,y\}$ is Leinert, as the only way for $x_{i}\neq x_{i+1}$ is if $x_i=x$ for all odd $i$ and $x_{i+1}=y$, and vice versa.
     \item Similarly to before, $G=\mathbf{F}_1\times\mathbf{F}_2$ is also Leinert.
 \end{enumerate}

 \subsection{Creating upper bounds}
 A primary method of showing that $X$ is statistically Leinert is to construct an upper bound on $b_{2n}$. Here is one way to count the number of bad sets, possibly construct upper bounds for $b_{2n}$ in the case $G=\mathbf{F}_{s_1}\times\mathbf{F}_{s_2}$:
 \begin{enumerate}
     \item Consider the string from the previous subsection,
     \[S=x_1^{-1}x_2y_1^{-1}y_2x_2^{-1}x_1y_2^{-1}y_1=e.\]
     We consider all strings of this length and form that equal the identity. We have to pick two different generators each (and we can swap the $x$'s and $y$'s in their role), which leads to a total of $2s_1(s_1-1)s_2(s_2-1)$ different strings of that form. Let's call this form the kernel string. 
     \item An alternative upper bound is the following, which is a very loose bound. The total number of walks of length $2n$ on the free group of $s$ generators from the identity to the identity is bounded by $(2s-1)^{n}\binom{2n}{n}$.
 \end{enumerate}
 Consider $s_1=s_2=2$. This does not form a Leinert set due to the arguments from before, 

\section{Numerical Algorithms}
 \subsection{Some numerics}
 For some numerics on the number of bad strings in various cases, consider the group $G$ from before. For strings of length $2n<8$, there are no bad strings. Here are some computations on the exact number of bad strings, or at least tight lower bounds. For notation purposes, if $S$ is a valid string of the form
 \[S=x_{1}^{-1}x_2\ldots x_{2n-1}^{-1}x_{2n},\]
 then define the conjugate string $\tilde S$ to be
 \[\tilde S=x_{1}x_2^{-1}\ldots x_{2n-1}x_{2n}^{-1}.\]
 \begin{itemize}
     \item For $2n=8$, we established that $b_{8}=2s_1(s_1-1)s_2(s_2-1)$.
     \item For $2n=10$, we can use every single string from before and append elements to the end that cancel each other out. Using the kernel form for $S_8$ (string of length 8), we can let $\tilde S = z_1^{-1}S_{8}z_1$, where $z_1\in X$ and $z_1\neq x_1$ or $y_2$. So we have $s-2$ choices for $z_1$. We can also consider strings of the form
     \[S_{10}=x_1^{-1}w_2x_4^{-1}y_1y_2^{-1}x_4w_2^{-1}x_1y_2^{-1}y_1=e.\]
     \item For $2n=12$, we can use any string from $b_{10}$ or use strings of the following form:
     \[S_{12}=x_1^{-1}w_2w_3^{-1}x_4y_1^{-1}y_2x_4^{-1}w_3w_2^{-1}x_1y_2^{-1}y_1=e.\]
     Given a kernel string $S_8$, we can look at strings of the form $z_1^{-1}z_2S_8z_2^{-1}z_1$, where we choose $z_1$ and $z_2$. $z_2$ has $s-2$ options, so $z_1$ has $s-3$ options. For strings of the form $S_{10}$, we have $b_8$ times the number of options for $w_2$ and $w_3$, which is hard to determine.
     \begin{itemize}
         \item If $w_2=x_4$, then $w_3$ has $s-1$ options. If $w_3=x_1$, then $w_2$ has $s-1$ options.
         \item Otherwise, $w_2$ has $s-2$ options from not equaling $x_4$ and $w_3$ has $s-2$ options as well.
     \end{itemize}
     We multiply by 2 again since we can also place the $w$'s between the $y$'s instead, so we have
     \[b_{12}=b_8[(s-2)(s-3) + 4(s-1) + 2(s-2)^2]\]
     \item In general, given $b_{2n}$, we have two options:
     \begin{enumerate}
         \item We can take a bad string $S_{2n-2}$ and conjugate by a letter $z$, giving $(s-2)$ strings.
         \item We can produce a bad string that has no bad substrings. The number of these is hard to find.
    \end{enumerate}
    
 \end{itemize}
\subsection{Growth of Bad Strings}
For the purposes of testing strings, a few tests were run. The numerical algorithm randomized a string of length 2m, using n generators. 

Here we shall consider five distinct cases: 
\begin{enumerate}
    \item 2 generators, 100 samples/length, length 1 to 40
    \item 2 generators, 1000 samples/length, length 1 to 40
\end{enumerate}

Three tests were applied: \begin{enumerate}
    \item Parity Test: The strings were then tested via a parity test to determine whether a necessary but not sufficient condition was met: there were an equal number of generators and their corresponding inverse elements.
    \item Adjacent Repeating Element Test: In this test adjacent elements in a string were checked to see if there were any repetitions, i.e. 'aa' or 'xx' etc.
    \item Reduction and Reordering Test: After the first two tests were completed, the strings were reordered cyclically and then reduced by checking whether adjacent elements were inverses of each other, in which case they were replaced by the identity element.
\end{enumerate}

In the first cycle of sample runs (not including the Adjacent Repeating Element Test), the following decay rate was observed.

\includegraphics[scale = .35]{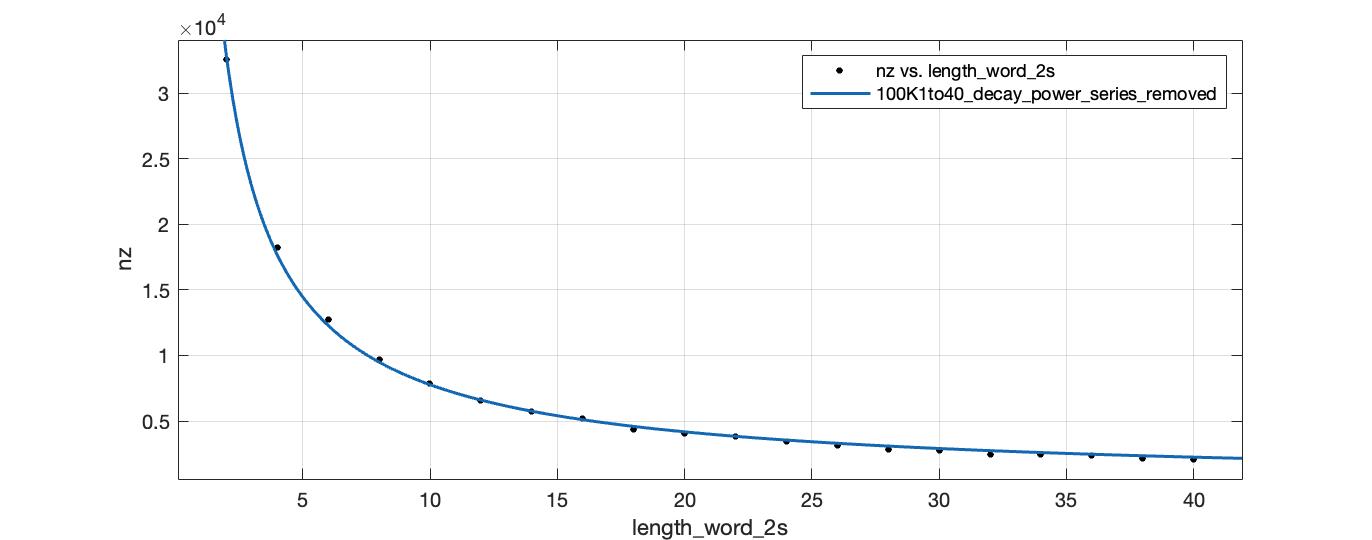}

In the second cycle of sample runs (including the Adjacent Repeating Element Test), the following decay rate was observed.

A Monte Carlo simulation was then applied to reduce computation time. The algorithm is presented in the appendix.
\section{Result and Numerical Verification}
When studying the radius of convergence of product of free groups, a specific instance of almost Leinert set.  We found a function Q(z, G(z)) that incorporates the number of bad strings in the set such that:
\begin{equation}
    P(zG(z)) \le G(z) \le Q(z, G(z))
\end{equation}
Here \begin{equation}
    P(t) = 1+\frac{1}{2}\sum_{x \in X}(\sqrt{1+4|\alpha(x)|^2t^2}-1)
\end{equation} and \begin{equation}
    Q(t,z) = 1 + \frac{1}{2}\sum_{i,j}(\sqrt{(1+\frac{zt^2}{c^2-t^2})^2+4\alpha_{i,j}^2z^2t^2}+\frac{zt^2}{c^2-t^2}-1)
\end{equation}
In our numerical calculations with increasing dimensions in the group formed by the product of two free groups with two generators, the lower bound on radius of convergence for G(z) differed from the upper bound by ~3. 
\begin{center}
    \includegraphics[width=0.45\textwidth]{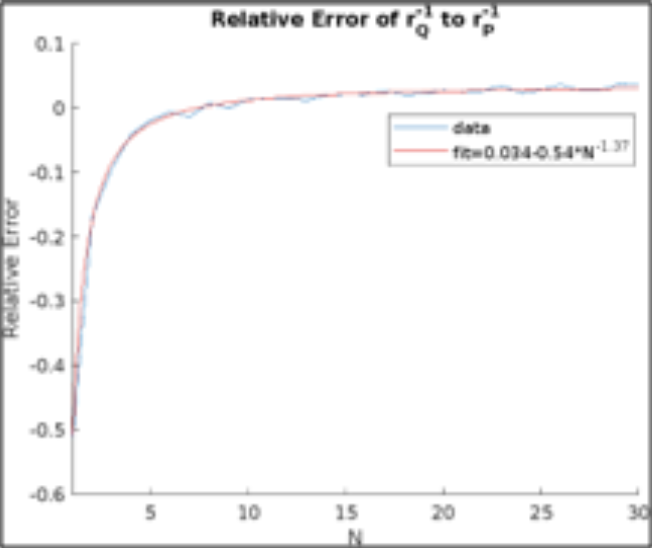}
\end{center}
\section{Conclusion and Future Works}
The results seem to be following our intuition created by past studies. This new function has been successful in tightly bounding the radius of convergence of almost Leinert Sets with the radius of convergence found by Woess for Leinert Sets.

Future works should consider analytically exploring a tighter bound on the "irreducible" strings as the size of the string exponentially decrease in proportion to valid strings. For the purposes of this project, only bad strings up to length sixteen were found, and the minimum length of the string found was length 8. As the length increases and bad strings are added in, the number of bad strings are seen to exponentially increase in proportion to valid strings.

\section{Acknowledgements}
Special acknowledgements to Dr. Thomas J. Sinclair, Purdue University Faculty, for mentoring this project.
This work was partially funded by NSF grant DMS-2055155.

\section{Appendix: Algorithm}

The algorithm mentioned in section 6.2 is presented here:
\begin{algorithmic}
    \State \textbf{Algorithm:} badStringCalculator 
    \State Define alpha be generator set
    \State Define the length of each word length\_word
    \State Define the number of words num\_words
    \State Define iter\_words as a vector 
    \State Define num\_reduced as a matrix 1 by iter\_words 
    \State Define num\_bad as a matrix 1 by iter\_words 
    \State Instantiate the counter number of reduced words num\_red\_words and the counter for the number of bad strings num\_bad\_strings
    \State Define inst 
    \For {i in [1, length(iter\_words)]}
        \For {j in [1:iter\_words(i)]}
        \State Code to reduce and reorder each word by verifying consecutive letters are inverses of each other
        \State word split 
        \State Define word as reducing placeholder red\_word, and character place holder temp\_char, and a steps counter steps
        \State Define matrix to compare consecutive letters parity\_mat
        \State new\_word = word;
        \State insert 

        \For{k in range [1, length(word)]}
            \State Work by case each possible letter case to fill the parity matrix, x is a +1 in the first column while y (its inverse) is -1 in the first column. Similarly, a and b are +1 and -1, respectively, in the second column.
            \If{word(k) == "a"}
                parity\_mat(1) = parity\_mat(1)+1;
            \ElsIf{word(k) == "b"}
                parity\_mat(1) = parity\_mat(1)-1;
            \ElsIf{word(k) == "x"}
                parity\_mat(2) = parity\_mat(2)+1;
            \ElsIf{word(k) == "y"}
                parity\_mat(2) = parity\_mat(2)-1;
            \EndIf
        \EndFor
        \State Define nn as an empty string
        \State Define vector of zeros the same length as word called nn1
        \For{k = 1:length(word)}
        \State Create a copy of word as nn and use nn1 to turn word into an array of letters, essentially splitting the string.
            \If{word(k) == "a"}
                nn = append(nn,"a");
                nn1(k) = 'a';
            \ElsIf{word(k) == "b"}
                nn = append(nn,"b");
                nn1(k) = 'b';
            \ElsIf{word(k) == "x"}
                nn = append(nn,"x");
                nn1(k) = 'x';
            \ElsIf{word(k) == "y"}
                nn = append(nn,"y");
                nn1(k) = 'y';
            \EndIf
        \EndFor
        \State Create reduced word to be a copy of nn. This copy will be parsed to reduce if inverses are adjacent to one another.
        \State Define steps to be an integer which counts how many reductions have been done.
        \State Define num\_red\_word to be an empty array.
        \State Define num\_bad to be an empty array.
        \If{parity\_mat == zeros(size(parity\_mat))}
            \While{red\_word ~= "" \&\& steps <= 10}
                \If{contains(red\_word,'ab')}
                    red\_word = replace(red\_word,'ab',' ')
                    inst = count(red\_word,' ');
                    red\_word = strrep(red\_word, ' ', '');
                    steps = steps+1;
                \EndIf
                \If{contains(red\_word,'ba')}
                    red\_word = replace(red\_word,'ba',' ')
                    inst = count(red\_word,' ');
                    red\_word = strrep(red\_word, ' ', '');
                    steps = steps+1;
                \EndIf
                \If{contains(red\_word,'xy')}
                    red\_word = replace(red\_word,'xy',' ')
                    inst = count(red\_word,' ');
                    red\_word = strrep(red\_word, ' ', '');
                    steps = steps+1;
                \EndIf
                \If{contains(red\_word,'yx')}
                    red\_word = replace(red\_word,'yx',' ')
                    inst = count(red\_word,' ');
                    red\_word = strrep(red\_word, ' ', '');
                    steps = steps+1;
                \EndIf
            \EndWhile
        \ElsIf{parity\_mat ~= zeros(size(parity\_mat))}
            num\_bad\_strings = num\_bad\_strings+1;
        \EndIf 
        num\_reduced(i) = num\_red\_word;
        num\_bad(i) = num\_bad\_strings
        \EndFor 
    \EndFor 
\State Define samples to be run. Instantiated to 1000.
\State Define length as an int variable for string.
\State Now we need to define the distribution to sample from for the Monte-Carlo sampling process.
\State Define mu to be 0.23 as the mean of the normal distribution.
\State Define sigma to be 0.25 as the standard deviation of the normal distribution.
\State Define L to be the length of iter\_words
\State Define Rmat to be a matrix of NaN values of size L x samples.
\State Define pdfNormalCd to be the normal probability distribution function of the normal distribution with mean mu and standard deviation sigma, evaluated at the values in iter\_words.
\For{i = 1:samples}
    R=normrnd(mu,sigma,[L,1]);
    Rmat(:,i)=R;
\EndFor
\State Plot iter\_words vs num\_reduced and Plot iter\_words vs num\_bad on the same plot
\end{algorithmic}
\bibliographystyle{plainnat}
\bibliography{USE_OF_STATISTICALLY_LEINERT_SETS_TO_CALCULATE_RETURN_PROBABILITIES_OF_RANDOM_WALKS} 

\begin{thebibliography}{5}
\providecommand{\natexlab}[1]{#1}
\providecommand{\url}[1]{\texttt{#1}}
\expandafter\ifx\csname urlstyle\endcsname\relax
  \providecommand{\doi}[1]{doi: #1}\else
  \providecommand{\doi}{doi: \begingroup \urlstyle{rm}\Url}\fi

\bibitem[Akemann and Ostrand(1976)]{akemann}
Charles~A. Akemann and Phillip~A. Ostrand.
\newblock Computing norms in group c*-algebras.
\newblock \emph{American Journal of Mathematics}, 98\penalty0 (4):\penalty0
  1015--1047, 1976.
\newblock ISSN 00029327, 10806377.
\newblock URL \url{http://www.jstor.org/stable/2374039}.

\bibitem[Hastings(2007)]{hastings}
M.~B. Hastings.
\newblock Random unitaries give quantum expanders.
\newblock \emph{Phys. Rev. A}, 76:\penalty0 032315, Sep 2007.
\newblock \doi{10.1103/PhysRevA.76.032315}.
\newblock URL \url{https://link.aps.org/doi/10.1103/PhysRevA.76.032315}.

\bibitem[Mathworld(2023)]{cayley}
Wolfram Mathworld.
\newblock Cayley tree, 2023.
\newblock URL \url{https://mathworld.wolfram.com/CayleyTree.html}.
\newblock Accessed on May 18, 2023.

\bibitem[Woess(1986)]{woess}
Wolfgang Woess.
\newblock A short computation of the norms of free convolution operators.
\newblock \emph{Proceedings of the American Mathematical Society}, 96\penalty0
  (1):\penalty0 167--170, 1986.

\bibitem[Zummo(2022)]{uchicago}
J.~Zummo.
\newblock An introduction to geometric group theory.
\newblock \emph{University of Chicago REU}, Aug 2022.
\newblock URL \url{http://math.uchicago.edu/~may/REU2022/REUPapers/Zummo.pdf}.

\end{thebibliography}

\end{document}